\documentclass[10pt,a4paper]{article}
\usepackage{times,latexsym,amsfonts}

\def\bmatrix#1{\left[\matrix{#1}\right]}

\newcommand{\QED}{\hspace*{\fill}\mbox{\footnotesize$\Box$}\par\bigskip}
\newcommand{\QEDBBOX}{\hspace*{\fill}\rule{0.4em}{0.4em}}

\def\kmtheoremdef#1#2#3#4{
  \newtheorem{#1x}{#2}[section]
  \newenvironment{#1}{\begin{#1x}#3\relax}{#4\end{#1x}}}

\kmtheoremdef{assumption}{Assumption}{\rm}{}
\kmtheoremdef{corollary}{Corollary}{\it}{}
\kmtheoremdef{definition}{Definition}{\rm}{}
\kmtheoremdef{example}{Example}{\rm}{}
\kmtheoremdef{lemma}{Lemma}{\it}{}
\kmtheoremdef{note}{Note}{\rm}{}
\kmtheoremdef{proposition}{Proposition}{\it}{}
\kmtheoremdef{theorem}{Theorem}{\it}{}

\newenvironment{proof}{\par\noindent{\textbf{Proof}}.\ }{\QED}
\newenvironment{proofof}[1]{\par\noindent{\bf Proof of #1}.\ }{\QED}

\newcommand{\IR}{\mathbb{R}}
\newcommand{\Z}{\mathbb{Z}}
\newcommand{\Q}{\mathbb{Q}}

\begin{document}
\begin{center}
\textbf{\Large
  Feedback Stabilization over Commutative Rings\\
  with no Right-/Left-Coprime Factorizations
}

\bigskip
\bigskip
\begin{large}
  Kazuyoshi MORI$\dag\ddag$
\end{large}

$\dag$ Department of Electrical Engineering, Faculty of Engineering,\\
       Tohoku University, Sendai 980-8579, JAPAN\\
       (\texttt{Kazuyoshi.MORI@IEEE.ORG})
\end{center}

\renewcommand{\thefootnote}{\fnsymbol{footnote}} \footnotetext[3]{
The author wishes to thank to Dr.~A.~Quadrat (Centre d'Enseignement et
de Recherche en Math\'{e}matiques, Informatique et Calcul
Scientifique, ENPC, France) for valuable comments.  A~significant part
of this work was done during the author's stay at Institut de
Recherche en Cybern\'etique de Nantes in France, to which the author
is grateful for hospitality and support.}
\renewcommand{\thefootnote}{\arabic{footnote}}

\section{Introduction}
Anantharam showed in~\cite{bib:anantharam85a} the existence of a~model
in which some stabilizable plants do not have its right-/left-coprime
factorizations.  In this paper, we give a~condition of the
nonexistence of the right-/left-coprime factorizations of stabilizable
plants as a~generalization of Anantharam's result.  As examples of the
models which satisfy the condition, we present two models; one is
Anantharam's example and the other the discrete finite-time delay
system which does not have the unit delay.  We illustrate the
construction of stabilizing controllers of stabilizable single-input
single-output plants of such models.  The method presented here is an
application of the result of the necessary and sufficient condition of
the stabilizability over commutative rings, which has recently been
developed by Abe and the author\cite{bib:mori98bs2}.

\section{Preliminaries}
\label{S:Preliminary}
The reader is referred to Section\,2 of~\cite{bib:mori99ata} for the
notations of commutative rings, matrices, and modules commonly used
throughout the paper, for the formulation of the feedback
stabilization problem, and for the related previous results.

\section{Anantharam's Result and Its Generalization}
\subsection{Anantharam's Result}
In~\cite{bib:anantharam85a}, Anantharam considered the case where
$\Z[\sqrt{5}i]$ ($\simeq\Z[x]/(x^2+5)$) is the set of the stable causal
transfer functions, where $\Z$ is the ring of integers and~$i$ the
imaginary unit; that is, ${\cal A}=\Z[\sqrt{5}i]$.  The set of all possible
transfer functions is given as the field of fractions of~${\cal A}$; that
is, ${\cal F}=\Q(\sqrt{5}i)$.  He considered the single-input single-output
case and showed that the plant $p=(1+\sqrt{5}i)/2$ does not have its
coprime factorization over~${\cal A}$ but is stabilizable.  As a~result, the
question ``\emph{Is it always necessary that the plant and its
stabilizing controller individually have coprime factorizations when
the closed loop is stable?}''  posed in~\cite{bib:vidyasagar82a} was
negatively answered over general commutative rings.

\subsection{A Generalization}
The result given by Anantharam\cite{bib:anantharam85a} can be
generalized over commutative rings.
\begin{proposition}\label{Prop:07.Mar.99.134855}
If there exist $a,b,a',b'\in{\cal A}$ satisfying the following three
statements, then there exists a~causal stabilizable plant~$P$ which
does not have right-/left-coprime factorizations:
\begin{itemize}
\item[(i)] The equality $ab=a'b'$ holds.
\item[(ii)] The ratio $a/a'$ is a~causal transfer
function and does not have coprime factorization.
\item[(iii)] The pair $(a,b)$ is coprime over~${\cal A}$.
\QEDBBOX
\end{itemize}
\end{proposition}
To prove the proposition above, we use the following proposition for
the single-input single-output case:
\begin{proposition}\label{Prop:06.Mar.99.215231}
Let $a,b,a',b'\in{\cal A}$.  Suppose that $ab=a'b'$ holds.  Suppose further
that $a/a'$ is a~causal plant~$p$.  Then, even if the plant~$p$ does
not have a~coprime factorization, if $(a,b)$ is coprime, then the
plant~$p$ is stabilizable.
\end{proposition}
\begin{proof}
We employ the notations in Definition\,2.4 of~\cite{bib:mori99ata}.
Thus, ${\cal I}=\{\{ 1\},\{ 2\}\}$, say $I_1=\{ 1\}$ and $I_2=\{ 2\}$.  Then the
generalized elementary factors of the plant are given as follows:
\begin{eqnarray}
  \Lambda_{pI_1}&=&\{\lambda\in{\cal A}\,|\,\exists K\in{\cal A}^{(m+n)\times m}\  \lambda T  =K \Delta_{I_1} T\}\nonumber\\
           &=&\{\lambda\in{\cal A}\,|\,\lambda a'a^{-1}\in{\cal A}\},\nonumber\\
  \Lambda_{pI_2}&=&\{\lambda\in{\cal A}\,|\,\exists K\in{\cal A}^{(m+n)\times m}\  \lambda T  =K \Delta_{I_2} T\}\nonumber\\
           &=&\{\lambda\in{\cal A}\,|\,\lambda aa'^{-1}\in{\cal A}\}.
                            \label{E:06.Mar.99.175111}
\end{eqnarray}
It is obvious that $a\in \Lambda_{pI_1}$.  On the other hand, $b$ is a~member
of $\Lambda_{pI_1}$ since the ratio $aa'^{-1}$ in
(\ref{E:06.Mar.99.175111}) can be rewritten as $b'b^{-1}$.  Since
$(a,b)$ is coprime, we have $\Lambda_{pI_1}+\Lambda_{pI_2}={\cal A}$.  Hence by
Theorem\,2.1 in~\cite{bib:mori99ata} the plant is stabilizable.
\end{proof}
\begin{proofof}{Proposition\,\ref{Prop:07.Mar.99.134855}}
From the conditions~(i) and~(iii), and
Proposition\,\ref{Prop:06.Mar.99.215231}, the plant $a/a'$ is
stabilizable and does not have right-/left-coprime factorization by
the condition\,(ii).
\end{proofof}
We note that Proposition\,\ref{Prop:07.Mar.99.134855} gives
a~condition of the \emph{nonexistence\/} of the doubly coprime
factorization of the stabilizable plant.  On the other hand, it should
be noted that Sule in~\cite{bib:sule94a} has given a~condition of the
existence of the doubly coprime factorization of the stabilizable
plant, which is expressed as follows:
\begin{proposition}\label{Prop:07.Mar.99.163819}
\textup{\rm (Theorem\,3 of~\cite{bib:sule94a})~} Let $\max {\cal A}$ be
Noetherian and $\dim\max{\cal A}=0$, where $\max{\cal A}$ denotes the set of all
maximal ideals of~${\cal A}$.  Then the plant is stabilizable if and only if
it has a~doubly coprime factorization.
\QEDBBOX
\end{proposition}
In particular, if the plant is of the single-input single-output and
if~${\cal A}$ is a~unique factorization domain, then Raman and Liu
in~\cite{bib:raman84a} gave the result that the plant is stabilizable
if and only if it has a~doubly coprime factorization.

We now give two examples which satisfy the condition of
Proposition\,\ref{Prop:07.Mar.99.134855}.
\begin{example}\label{Ex:06.Mar.99.223139:1}
In the example given in~\cite{bib:anantharam85a}, that is, the case
where $p=(1+\sqrt{5}i)/2$, the numbers $1+\sqrt{5}i$, $1-\sqrt{5}i$,
$2$, $3$ are corresponding to the variables~$a$, $b$, $a'$, $b'$,
respectively, in Proposition\,\ref{Prop:07.Mar.99.134855}.  \QEDBBOX
\end{example}
We can make analogous examples.  For example, let ${\cal A}
=\Z[\sqrt{xy-1}i]$, where the integers~$x$ and~$y$ satisfy the
following conditions:
(i)~$\gcd(x,y)=1$ over $\Z$,
(ii)~$y>x\geq 2$,
(iii)~$xy-1$ is not square.
If the plant $p=(1+\sqrt{xy-1}i)/x$, then the numbers
$1+\sqrt{xy-1}i$, $1-\sqrt{xy-1}i$, $x$ and~$y$ are corresponding
to~$a$, $b$, $a'$, and~$b'$, respectively, in
Proposition\,\ref{Prop:07.Mar.99.134855}.
\begin{example}\label{Ex:06.Mar.99.223139:2}
Let us consider the discrete finite-time delay system.  On some
high-speed electronic circuits such as computer memory devices, they
cannot often have nonzero small delays.  We suppose here that the
system cannot have the unit delay as a~nonzero small delay.  In this
case, the set~${\cal A}$ becomes the set of polynomials generated by~$x^2$
and~$x^3$, that is, ${\cal A}=\IR[x^2,x^3]$, where~$x$ denotes the unit
delay operator.  Then~${\cal A}$ is not a~unique factorization domain but
a~Noetherian domain.  The set~${\cal Z}$ used to define the causality is
given as the set of polynomials in $\IR[x^2,x^3]$ whose constant terms
are zero; that is, ${\cal Z}=\{\alpha x^2+\beta x^3\,|\,\alpha,\beta\in{\cal A}\}$.

Let us suppose that $p=(1-x^3)/(1-x^2)\in{\cal P}$%
.  Since
$(1-x^3)(1+x^3)=(1-x^2)(1+x^2+x^4)$, the plant can be also expressed
as $p=(1+x^2+x^4)/(1+x^3)$.  Then $(1-x^3)$, $(1+x^3)$, $(1-x^2)$, and
$(1+x^2+x^4)$ are corresponding to $a,b,a',b'$ in
Proposition\,\ref{Prop:06.Mar.99.215231}.  Hence the plant does not
have its coprime factorization but is stabilizable.  \QEDBBOX
\end{example}

\section{Stabilizability and Construction of Stabilizing Controllers}
In this section we present first the stabilizability for models of
Examples\,\ref{Ex:06.Mar.99.223139:1} and~\ref{Ex:06.Mar.99.223139:2},
and then the construction of stabilizing controllers plants.

\subsection{Stabilizability}
The following two propositions give the stabilizability of all
transfer functions in both cases of ${\cal A}=\Z[\sqrt{5}i]$ and ${\cal A}
=\IR[x^2,x^3]$.
\begin{proposition}\label{Prop:09.Mar.99.150521}
Let ${\cal A}=\Z[\sqrt{5}i]$.  Suppose that ${\cal Z}=\{ 0\}$, so that ${\cal P}={\cal F}$.
 Then any transfer functions in~${\cal F}$ are stabilizable.
\end{proposition}
\begin{proof}
In the case $p=0$, the plant~$p$ is obviously stabilizable.  Hence in
the following we assume without loss of generality that $p\neq 0$.  Let
us consider a~plant~$p$ is expressed as
  $p=(\alpha_1+\alpha_2\sqrt{5}i)/\beta$,  
where $\alpha_1,\alpha_2,\beta\in\Z$.  Let~$g$ denote $\gcd(\alpha_1^2+5 \alpha_2^2,\beta)$
over $\Z$.  Let~$\alpha'$ be an integer such that $\alpha_1^2+5 \alpha_2^2=\alpha'g$.
We note that~$\alpha'$ does not have~$\beta$ as a~factor.

To show the stabilizability, we here apply (iii) of Theorem\,2.1
in~\cite{bib:mori99ata}.  Since the plant is of the single-input
single-output, the set~${\cal I}$ defined in Definition\,2.4
of~\cite{bib:mori99ata} is equal to $\{\{ 1\},\{ 2\}\}$, say $I_1=\{ 1\}$
and $I_2=\{ 2\}$.  Then the generalized elementary factors $\Lambda_{pI_1}$
and $\Lambda_{pI_2}$ are given as follows:
\begin{equation}\label{E:09.Mar.99.150544}
  \Lambda_{pI_1} = \{\lambda\in{\cal A}\,|\,\lambda dn^{-1}\in{\cal A}\},
  \Lambda_{pI_2} = \{\lambda\in{\cal A}\,|\,\lambda nd^{-1}\in{\cal A}\}.
\end{equation}
Then one can check that $\alpha'\in \Lambda_{pI_1}$ and $\beta\in \Lambda_{pI_2}$.
Since $\gcd(\alpha',\beta)=1$ over $\Z$, we have $\Lambda_{pI_1}+\Lambda_{pI_2}={\cal A}
$.  Hence any plant~$p$ is stabilizable by Theorem\,2.1 in~\cite{bib:mori99ata}.
\end{proof}

\begin{proposition}\label{Prop:09.Mar.99.151101}
In the case of ${\cal A}=\IR[x^2,x^3]$, any causal transfer functions are
stabilizable.
\end{proposition}
\begin{proof}
Suppose that~$p$ is a~causal plant in~${\cal P}$ expressed as $p=n/d$ with
$n\in{\cal A}$ and $d\in{\cal A}\backslash{\cal Z}$.  Let~$g$ be $\gcd(n,d)$ over $\IR[x]$ rather
than over~${\cal A}$.  We assume without loss of generality that~$g$ is
expressed as $1+g_1x$ with $g_1\in\IR$%
.

Let~$n'$ and~$d'$ be polynomials $n/g$ and $d/g$, respectively, in
$\IR[x]$ and further~$n''$ and~$d''$ be polynomials in~${\cal A}$ defined as
follows:
\begin{eqnarray*}
  n''&=&n' \textstyle (1+g_1x+\frac{2}{9}g_1^2x^2)\\
     &=&\left\{
          \begin{array}{ll}
           n' \textstyle (1+\frac{1}{3}g_1x)(1+\frac{2}{3}g_1x) 
                              &{\ \ } \mbox{(if $g_1\neq 0$)},\\
           n                  &{\ \ } \mbox{(if $g_1=0$)},
          \end{array}
        \right.\\
  d''&=&d' \textstyle (1+g_1x+\frac{2}{9}g_1^2x^2)\\
     &=&\left\{
          \begin{array}{ll}
           d' \textstyle (1+\frac{1}{3}g_1x)(1+\frac{2}{3}g_1x) 
                              &{\ \ } \mbox{(if $g_1\neq 0$)},\\
           d                  &{\ \ } \mbox{(if $g_1=0$)}
          \end{array}
        \right.
\end{eqnarray*}
(Note here that $\frac{2}{9}$ above can be other values).  Then one
can check that $n\in \Lambda_{pI_1}$, $d''\in \Lambda_{pI_2}$ and the greatest
common divisor of~$n$ and~$d''$ over $\IR[x]$ is a~unit.  Let~$\alpha$, $
\beta$ be in $\IR[x]$ such that the following equation holds over
$\IR[x]$:
\begin{equation}\label{E:08.Mar.99.135111}
    (\alpha+rd'')n+(\beta-rn)d''=1,
\end{equation}
where~$r$ is an arbitrary element in $\IR[x]$.  Let~$n_0$, $d_0$, $\alpha_
0$, and~$\beta_0$ denote the constant terms of~$n$, $d''$, $\alpha$, and~$\beta
$, respectively.  Similarly let~$\alpha_1$, $\beta_1$, $\alpha_1'$, and~$\beta_1'$
denote the coefficients of~$\alpha$, $\beta$, $\alpha+rd''$, and $\beta-rn$,
respectively, with the degree~$1$ with respect to the variable~$x$.
In order to show that (\ref{E:08.Mar.99.135111}) holds over~${\cal A}$, we
want to find~$r$ such that $\alpha_ 1'=\beta_1'=0$.  To do so, we let $r=(\alpha_
0 \beta_1-\alpha_1 \beta_ 0)x$.  Then it is easy to check that $\alpha_ 1'=\beta_1'=0$
hold from the relations that $\alpha_0n_0+\beta_0d_0=1$ and $\alpha_1n_0+\beta_
1d_0=0$.  Now that $ (\alpha+rd'')$, $(\beta-rn)\in{\cal A}$, we have $\Lambda_{pI_1}+\Lambda_
{pI_2}={\cal A}$, so that the plant is stabilizable.
\end{proof}

\subsection{Construction of Stabilizing Controllers}
We present here the method to construct stabilizing controllers under
the cases (i) ${\cal A}=\Z[\sqrt{5}i]$ or $=\IR[x^2,x^3]$ and (ii)
single-input single-output plant.  This is an application of the proof
(``(iii)$\rightarrow$(i)'' part) of Theorem\,2.1 in~\cite{bib:mori99ata} (for the proof,
see~\cite{bib:mori98bs2}).

Since the plant~$p$ is of the single-input single-output, the set~${\cal I}$ defined
in Definition\,2.4 of~\cite{bib:mori99ata} is equal to $\{\{ 1\},\{ 2\}\}
$, say $I_1=\{ 1\}$ and $I_2=\{ 2\}$ as in the proof of
Proposition\,\ref{Prop:09.Mar.99.150521}.  Two generalized elementary
factors $\Lambda_{pI_1}$ and $\Lambda_{pI_2}$ are given as
(\ref{E:09.Mar.99.150544}).  Since in the cases ${\cal A}=\Z[\sqrt{5}i]$ and
${\cal A}=\IR[x^2,x^3]$ any causal transfer functions are stabilizable, $\Lambda_
{pI_1}+\Lambda_{pI_2}={\cal A}$ holds by Theorem\,2.1 in~\cite{bib:mori99ata}.
We should find $\lambda_{I_1}\in \Lambda_{pI_1}$ and $\lambda_ {I_2}\in \Lambda_{pI_2}$ such
that $\lambda_{I_1}+\lambda_ {I_2}=1$.

From Lemmas\,4.7 and 4.10 of~\cite{bib:mori98bs2}, there exist
(right-/left-)coprime factorizations over ${\cal A}_{\lambda_{I}}$ for $I=I_1$,
$I_2$.  We let $n_1=1$, $d_1=1/p$, $n_2=p$, $d_2=1$ with
$p=n_1/d_1=n_2/d_2$, $n_1,d_1\in{\cal A}_{\lambda_{I_1}}$, and $n_2,d_2\in{\cal A}_{\lambda_
{I_2}}$.  Then the coprime factorizations are obtained as follows:
\begin{eqnarray*}
&&  y_1n_1+x_1d_1=1 {\ \ } \mbox{(over ${\cal A}_{\lambda_{I_1}}$)}{\ \ } \mbox{and}\\
&&  y_2n_2+x_2d_2=1 {\ \ } \mbox{(over ${\cal A}_{\lambda_{I_2}}$)},
\end{eqnarray*}
where $y_1=1$, $x_1=0$, $y_2=0$ and $x_2=1$.  When we use the
parameters $r_1\in{\cal A}_{\lambda_{I_1}}$ and $r_2\in{\cal A}_{\lambda_ {I_2}}$, we also
have
\begin{eqnarray*}
&& (y_1+r_1d_1)\cdot n_1+(x_1-r_1n_1)\cdot d_1=1 {\ \ }\mbox{(over ${\cal A}_{\lambda_{I_1}}$)},~~~~~~\\
&& (y_2+r_2d_2)\cdot n_2+(x_2-r_2n_2)\cdot d_2=1 {\ \ }\mbox{(over ${\cal A}_{\lambda_{I_2}}$)}.~~~~~~
\end{eqnarray*}
From the proof of Theorem\,2.1 in~\cite{bib:mori99ata} (Theorem\,3.2
of \cite{bib:mori98bs2}), a~stabilizing controllers of the plant~$p$ is
given as the following form
\begin{equation}\label{E:06.Mar.99.234746}
  c=\frac{a_1 \lambda_{I_1}^{\omega}d_1(y_1+r_1d_1)+a_2 \lambda_{I_2}^{\omega}d_2(y_2+r_2d_2)}%
         {a_1 \lambda_{I_1}^{\omega}d_1(x_1-r_1n_1)+a_2 \lambda_{I_2}^{\omega}d_2(x_2-r_2n_2)},
\end{equation}
where~$\omega$ is a~sufficiently large positive integer and
$a_1,a_2,r_1,r_2$ are elements in~${\cal A}$ such that the following three
conditions hold:
\begin{itemize}
\item[(i)]  The equality $a_1 \lambda_{I_1}^{\omega}+a_2 \lambda_{I_2}^{\omega}=1$ holds.
\item[(ii)] The following matrices are over~${\cal A}$ for $k=1,2$:
\begin{equation}\label{E:01.Feb.99.102744}
  \begin{array}{l}
  a_k \lambda_{I_k}^{\omega} n_k(x_k-r_kn_k),~
  a_k \lambda_{I_k}^{\omega} n_k(y_k+r_kd_k),\hspace*{3em}\\
  \hspace*{3em}
  a_k \lambda_{I_k}^{\omega} d_k(x_k-r_kn_k),~
  a_k \lambda_{I_k}^{\omega} d_k(y_k+r_kd_k).
  \end{array}
\end{equation}
\item[(iii)] The denominator of~$c$ is nonzero.
\end{itemize}

In the following we continue Examples\,\ref{Ex:06.Mar.99.223139:1}
and~\ref{Ex:06.Mar.99.223139:2} in which the stabilizing controllers
are constructed.

\par\noindent
\textbf{Example \ref{Ex:06.Mar.99.223139:1} (continued)~} 
Let us stabilize the plant $p=(1+\sqrt{5}i)/2$.  We can find that 
   $\alpha'=3=:\lambda_{I_1}\in \Lambda_{pI_1}$ and $\beta=2=:\lambda_{I_2}\in \Lambda_{pI_2}$, 
where~$\alpha'$ and~$\beta$ are the symbols used in the proof of
Proposition\,\ref{Prop:09.Mar.99.150521}.  By using $\lambda_{I_1}$ and $\lambda_
{I_2}$ the stabilizing controllers are obtained in the form
(\ref{E:06.Mar.99.234746}).  For example, let us choose $r_1=r_2=0$.
Then we can select $\omega=1$ to satisfy the condition~(ii) above.  The
coefficients~$a_1$ and~$a_2$ are~$1$ and $-1$, respectively.  We
obtain a~stabilizing controller $c=\frac{-1+\sqrt{5}i}{2}$, which is
same as the stabilizing controller given in~\cite{bib:anantharam85a}.

Let us present the parameterization of the stabilizing controllers
according to~\cite{bib:mori99ata}.
Now $H_0=H(p,c)$ is expressed as
\[
  H_0=\bmatrix{ -2           & 1+\sqrt{5}i \cr
                1- \sqrt{5}i & -2 }.
\]
According to Theorem\,4.3 of~\cite{bib:mori99ata}, the set of all
$H(p,c)$'s with all stabilizing controllers $c$'s, denoted by ${\cal H}(p;{\cal A})
$, is given as follows:
\begin{eqnarray}
\lefteqn{\hspace*{1.0em}\mbox{\normalsize${\cal H}(p;{\cal A})=$}}\label{E:23.Apr.99.171348}\\
&&
   \left\{
   \begin{array}{l}
     \bmatrix{ h_{11} & h_{12}\cr
               h_{21} & h_{22}} \Biggm|
         q_{11},q_{12},q_{21},q_{22}\in{\cal A},\\
   
      ~~~h_{11}=h_{22}=
      3\sqrt{5}iq_{12}-2\sqrt{5}iq_{21}+6q_{11}~~\\
         \multicolumn{1}{r}{-3q_{12}-2q_{21}+6q_{22}-2,}\\
      ~~~h_{12}=-3\sqrt{5}iq_{11}+2\sqrt{5}iq_{21}-3\sqrt{5}iq_{22}~~\\
         \multicolumn{1}{r}{+\sqrt{5}i-3q_{11}+9q_{12}-4q_{21}-3q_{22}+1,}\\
      ~~~h_{21}=2\sqrt{5}iq_{11}-2\sqrt{5}iq_{12}+2\sqrt{5}iq_{22}-\sqrt{5}i~~\\
         \multicolumn{1}{r}{-2q_{11}-4q_{12}+4q_{21}-2q_{22}+1,}\\
      ~~~\mbox{ $h_{11}$ and $h_{22}$ are nonzero.}
   \end{array}
  \right\}.\nonumber
\end{eqnarray}
In the set ${\cal H}(p;{\cal A})$ in (\ref{E:23.Apr.99.171348}), $q_{ij}$ is
corresponding to the $(i,j)$-entry of the matrix~$Q$ used in~\S\,4
of~\cite{bib:mori99ata}.  One can observe that
$(1+\sqrt{5}i)h_{11}=-2h_{12}$, so that $p=-h_{12}h_{11}^{-1}$.  By
Corollary\,4.1 of~\cite{bib:mori99ata} any stabilizing controllers are
expressed as $h_{11}^{-1}h_{21}$ or equivalently $h_{21}h_{22}^{-1}$
provided that $h_{11}$ and $h_{22}$ are nonzero.

\par\noindent 
\textbf{Example \ref{Ex:06.Mar.99.223139:2} (continued)~} Let us
construct stabilizing controllers of the plant $p=(1-x^3)/(1-x^2)$.
In this case, $g$, $n$ and~$d''$ used in the proof of
Proposition\,\ref{Prop:09.Mar.99.151101} are
  $1-x$,
  $1-x^3$ and
  $1- \frac{7}{9}x^2+\frac{2}{9}x^3$,
respectively.  We calculate~$\alpha$ and~$\beta$ used in the proof satisfying
(\ref{E:08.Mar.99.135111}).  They are computed as
  $\alpha=- \frac{101}{988} - \frac{441}{988}x  +\frac{77}{494}x^2$ and
  $\beta= \frac{1089}{988} + \frac{441}{988}x +\frac{693}{988}x^2$.
By letting $r=\frac{441}{988}x$, all of $(\alpha+rd'')$, $(\beta-rn)$, $n$, $d''$ in
(\ref{E:08.Mar.99.135111}) become in~${\cal A}$.

Let $\lambda_{I_1}=n\in \Lambda_{pI_1}$ and $\lambda_{I_2}=d''\in \Lambda_{pI_2}$.  Then a~stabilizing controllers is given in the form of
(\ref{E:06.Mar.99.234746}).  For example, let us choose $r_1=r_2=0$.
Then we can select $\omega=1$ to satisfy the condition~(ii).  The
coefficients~$a_1$ and~$a_2$ are
      $\textstyle - \frac{101}{988} + \frac{77}{494}x^2 
                  - \frac{343}{988}x^3  +\frac{49}{494}x^4$, and 
      $\textstyle   \frac{1089}{988} + \frac{693}{988}x^2 
                  + \frac{441}{988}x^4$,
respectively.
Finally, from the formula of (\ref{E:06.Mar.99.234746}) we obtain
a~stabilizing controller given as follows:
\begin{small}
\[
c=\frac{
         - 101
         + 255x^2 
         - 343x^3 
         -  56x^4 
         + 343x^5 
         -  98x^6 }%
     {     1089
         - 154x^2 
         + 242x^3 
         -  98x^4 
         + 154x^5 
         - 343x^6 
         +  98x^7 }.
\]
\end{small}%
Again, let us consider the parameterization of the stabilizing
controllers.  Although we can state the set ${\cal H}(p;{\cal A})$, we do not
state it unlike (\ref{E:23.Apr.99.171348}) because of space
limitations.  Further any stabilizing controllers are expressed as
$h_{11}^{-1}h_{21}$ or equivalently $h_{21}h_{22}^{-1}$ provided that
$h_{11}$ and $h_{22}$ are nonzero.

\section{Concluding Remarks}
In this paper, we have presented a~generalization of Anantharam's
result and given a~condition of the nonexistence of the
right-/left-coprime factorizations of stabilizable plants.  As
examples satisfying the obtained condition, two models were presented.
We have also presented a~method to construct stabilizing controllers
of stabilizable single-input single-output plants of such models. 

Propositions\,\ref{Prop:07.Mar.99.134855} and
\ref{Prop:07.Mar.99.163819} give the conditions of the
existence/nonexistence of the doubly coprime factorizations of
stabilizable plants.  However they do not characterize the commutative
ring~${\cal A}$ on which there exist the doubly coprime factorizations of
stabilizable plants.  This problem is not solved yet.

%
%


\end{document}